\documentclass{amsart}
\usepackage[english]{babel}
\usepackage[latin1]{inputenc}
\usepackage{mathrsfs}
\usepackage{amsmath}
\usepackage{amsthm}
\usepackage{amssymb}
\usepackage{indentfirst}

\usepackage[dvips]{graphicx}

\usepackage{ifthen}
\newcommand{\cit}[1]{{\rm \textbf{#1}}}

\newcommand{\Ref}[2]{\cit{%
\ifthenelse{\equal{#1}{thm}}{Theorem}{}%
\ifthenelse{\equal{#1}{prop}}{Proposition}{}%
\ifthenelse{\equal{#1}{lem}}{Lemma}{}%
\ifthenelse{\equal{#1}{cor}}{Corollary}{}%
\ifthenelse{\equal{#1}{defn}}{Definition}{}%
\ifthenelse{\equal{#1}{oss}}{Remark}{}%
\ifthenelse{\equal{#1}{sec}}{Section}{}%
\ifthenelse{\equal{#1}{ex}}{Example}{}%
\ifthenelse{\equal{#1}{conj}}{Conjecture}{}%
\ifthenelse{\equal{#1}{ssec}}{Subsection}{}%
\ifthenelse{\equal{#1}{tab}}{Table}{}%
\ifthenelse{\equal{#1}{cla}}{Claim}{}%
\  \ref{#1:#2}%
}}

\theoremstyle{plain} 
\newtheorem{prop}{Proposition}[section]
\newtheorem{thm}[prop]{Theorem}

\newtheorem*{thmnn}{Theorem}
\newtheorem{lem}[prop]{Lemma} 
\newtheorem{cor}[prop]{Corollary}
\newtheorem{conj}[prop]{Conjecture}
\newtheorem*{conjnn}{Conjecture}

\theoremstyle{remark}

\newtheorem{oss}[prop]{Remark}
\newtheorem{ex}[prop]{Example}

\theoremstyle{definition}
\newtheorem{defn}[prop]{Definition}
\newcommand{\hk}{hyperk\"{a}hler }

\newcommand{\kahl}{K\"{a}hler }
\newcommand{\ktiposp}{$K3^{[2]}$ type }
\newcommand{\ktipo}{$K3^{[2]}$ type}
\newcommand{\kntipo}{$K3^{[n]}$ type}
\newcommand{\kntiposp}{$K3^{[n]}$ type }
\newcommand{\ie}{i.~e.~}

\begin{document}
\title[Symplectic automorphisms on manifolds of \kntipo]{Towards a classification of symplectic automorphisms on manifolds of $K3^{[n]}$ type}
\author{Giovanni Mongardi},
\address{Universit\`{a} degli studi di Milano, Dipartimento di Matematica, via Cesare Saldini, 50}
\email{giovanni.mongardi@unimi.it}
\begin{abstract}

The present paper is devoted to the classification of symplectic automorphisms of some \hk manifolds. The results contained here are an explicit classification of prime order automorphisms on manifolds of \kntiposp and a proof that all finite groups of symplectic automorphisms of such manifolds are contained in Conway's group $Co_1$.  

\end{abstract}
\thanks{Supported by FIRB 2012 ``Spazi di Moduli e applicazioni'', by SFB/TR 45 ``Periods, moduli spaces and arithmetic of algebraic varieties'' and partially by the Max Planck Institute in Mathematics.}

\maketitle

\section*{Introduction}
Automorphisms of symplectic manifolds have attracted a lot of interest, mainly in the two dimensional setting. Automorphisms of $K3$ surfaces preserving the symplectic form have been fully classified, starting from the foundational work of Nikulin \cite{nik1} and Mukai \cite{muk} and of several other eventually leads to a list of groups and actions on cohomology by Hashimoto \cite{hashi}.\\
In higher dimensions, far less is known: there are results concerning deformations of automorphisms, such as \cite{boi} by Boissi\`{e}re and \cite{me_natural}. There are also previous results in the peculiar case of involutions on four dimensional \hk manifolds by Camere \cite{cam}, O'Grady \cite{ogr2} and there is also \cite{me_invol}. Symplectic automorphisms have already been addressed in \cite{me_autom1} and the recent \cite{bcs} by Boissi\`{e}re, Camere and Sarti addresses automorphisms that do not preserve the symplectic form.\\ 
In the present paper, we specialise to manifolds obtained as smooth deformations of Hilbert schemes of $n$ points on $K3$ surfaces (shorthand: manifolds of \kntipo). Our results, obtained with techniques similar to Kondo's \cite{kon} approach for $K3$ surfaces, yield a classification of prime order automorphisms preserving the symplectic form and the following two characterisations of finite groups of symplectic automorphisms:
\begin{thmnn}
Let $X$ be a \hk manifold of \kntiposp and let $G$ be a finite group of symplectic automorphisms of $X$. Then $G\subset Co_1$ and $S_G(X)=S_G(\Lambda)$ for some conjugacy class of $G$ in $Co_1$. 
\end{thmnn}
See \Ref{thm}{sporadic} for details, here $\Lambda$ is the Leech lattice, \ie the unique even unimodular negative definite lattice with no elements of square $-2$ and $Aut(\Lambda)/\pm1=Co_1$. This can be seen as an analogous of Mukai's \cite{muk} result on $K3$ surfaces, where symplectic automorphisms were realised as subgroups of the Mathieu group $M_{23}\subset Co_1$. The Conway group appears also in some related fields, like derived autoequivalencies of $K3$ surfaces \cite[Huybrechts]{huy_der} or sigma models on $K3$ surfaces \cite[Gaberdiel, Hohenegger and Volpato]{ghv}. The following is a somewhat more technical sufficient condition:
\begin{thmnn}
Let $G\subset Co_0$ be a group of isometries of the Leech lattice with invariant sublattice $T$ of rank at least $4$ and the discriminant group of $T$ has less generators than its rank. Then there exist an integer $n$ and a manifold $X$ of \kntiposp such that $G\subset Aut_s(X)$.
\end{thmnn}
See \Ref{thm}{da_conway_ad_auto}. 
Based on the behaviour of some isometries not satisfying the above proposition, we conjecture the following.
\begin{conjnn}
There is a bijective correspondence between finite groups of symplectic automorphisms of manifolds of \kntiposp (for some $n$) and subgroups $G$ of $Co_1$ satisfying the conditions of the above proposition.
\end{conjnn}
Recently, Huybrechts \cite{huy_der} found a similar result for automorphisms coming from symplectic derived autoequivalencies of $K3$ surfaces. In \Ref{lem}{cond_huy} we prove that the sufficient condition he imposes on the group action is equivalent to ours.\\ 
The structure of the paper is as follows. In \Ref{sec}{k3n}, we gather several known results on manifolds of \kntipo. In \Ref{sec}{lattice} we collect some results about lattices, mainly about their discriminant groups. In \Ref{sec}{aut_in_co} we prove basic properties of invariant and coinvariant lattices associated to symplectic automorphisms and we use them to prove that finite groups of symplectic automorphisms can be embedded in the Conway group. The same result, but only in dimension four, was proven in \cite{me_autom1}. In \Ref{sec}{co_in_aut} we address the converse, namely how to obtain automorphisms from subgroups of the Conway group. We provide a sufficient and a necessary condition.
\Ref{sec}{prime} contains a classification of all symplectic automorphisms of prime order, up to their coinvariant lattice. In the special case of dimension four, we also classify the number of deformation classes of pairs $(X,\varphi)$ and the fixed locus of $\varphi$. Here $\varphi$ is a symplectic automorphism of prime order on a fourfold $X$ of \ktipo. In the case of an involution, this result is already proven in \cite{me_invol}.
Finally, in the appendix we state a series of known results about sublattices of the Leech lattice, which we will use in the previous sections. The results of the appendix are taken from many sources, mainly \cite{atlas}, \cite{con}, \cite{cur}, \cite{ero}, \cite{me_tesi}, \cite{me_autom1} and \cite{nip}\\

\section*{Notations}
Let $L$ be a lattice and let $G\subset O(L)$ be a group of isometries. We denote by $T_G(L)$ the $G$-invariant sublattice and with $S_G(L):=T_G(L)^\perp$ the coinvariant lattice. For any lattice $L$, we denote by $L(n)$ the lattice with the same $\mathbb{Z}$ module structure of $L$, but with quadratic form multiplied by $n$. If $X$ is a manifold of \kntiposp and $G\subset Aut(X)$, we will denote with $T_G(X)$ and $S_G(X)$ the invariant and coinvariant parts of the induced action of $G$ on $H^2(X,\mathbb{Z})$. In the same fashion we will call discriminant group of $X$ the group $H^2(X,\mathbb{Z})^\vee/H^2(X,\mathbb{Z})$. We will denote with $Aut_s(X)$ the group of automorphisms of a manifold of \kntiposp that preserve the symplectic form, \ie symplectic automorphisms. Any isometry between $H^2(X,\mathbb{Z})$ and a lattice $L$ is called a marking. In the present paper, all abstract groups are finite. This, in general, is not the case for automorphism groups of manifold of \kntipo, but we will always consider finite subgroups of these groups.

\section*{Acknowledgements}
The present paper is an improvement of my PhD thesis \cite{me_tesi} and of \cite{me_autom1}, and I would like to thank my advisor, K.G. O'Grady, for his support. I am also grateful to S. Boissi\`{e}re and A. Sarti for useful discussions and to D. Huybrechts for the discussion about his work \cite{huy_der}. Part of this work was carried out during the Hausdorff junior trimester program in algebraic geometry, whose support I would like to acknowledge.

\section{Manifolds of \kntipo}\label{sec:k3n}
In this section we give a short introduction to manifolds of \kntipo. They are one of the known series of \hk manifolds, with whom they share several properties.
\begin{defn}
Let $S$ be a $K3$ surface and let $S^{[n]}$ denote its Hilbert scheme of points. Let $X$ be a smooth deformation of $S^{[n]}$. Then $X$ is called of \kntipo.
\end{defn}
We remark that these manifolds inherit a symplectic form from the $K3$ surface, as shown in \cite{beau3}.
We will mainly need three properties of these manifolds, which are the Beauville-Bogomolov form, the global Torelli theorem and the structure of their \kahl cone.

\begin{thm}
Let $X$ be a manifold of \kntipo. Then there exists a canonically defined pairing $(\,,\,)_X$ on $H^2(X,\mathbb{C})$, the Beauville-Bogomolov pairing, and a constant $c_X$ (the Fujiki constant) such that the following holds:
\begin{equation}
(\alpha,\alpha)_X^n=c_x\int_X\alpha^{2n}.
\end{equation}

Moreover $c_X$ and $(\,,\,)_X$ are deformation and birational invariants.
\end{thm}
With this form, $H^2(X,\mathbb{Z})$ is isometric to the following
\begin{equation}\label{latticeK3n}
L_n:=U\oplus U\oplus U\oplus E_8(-1)\oplus E_8(-1)\oplus (2-2n).
\end{equation}
Where $U$ is the hyperbolic lattice, $E_8(-1)$ is the unique unimodular even negative definite lattice of rank $8$, $(2-2n)$ is $(\mathbb{Z},q)$ with $q(1)=2-2n$ and $\oplus$ denotes orthogonal direct sum.
Any such isometry is called a marking and, in analogy to the case of $K3$ surfaces, there exists a non Hausdorff moduli space of marked manifolds.
\begin{defn}
Let $(X,\phi)$ be a marked manifold of \kntipo. Let $\mathcal{M}_n$ be the set $\{(X,\phi)\}/\sim$ of marked manifolds of \kntiposp where $(X,\phi)\sim (X',\phi')$ if and only if there exists an isomorphism $f\,:X\,\rightarrow\,X'$ such that $f^*=\phi^{-1}\circ \phi'$.
\end{defn}

\begin{defn}
We define the period domain $\Omega_n$ as
\begin{equation}
\Omega_n=\{x\in \mathbb{P}(L_n\otimes\mathbb{C})\,|\,(x,x)_{L_n}=0,\,(x+\overline{x},x+\overline{x})_{L_n}>0\}.
\end{equation}
\end{defn}
The period map $\mathcal{P}\,:\,\mathcal{M}_n\rightarrow \Omega_n$ sends a pair $(X,f)$ to the line spanned by $f(\sigma_X)$.
\begin{defn}
Let $S$ be a $K3$ surface. Let $L_M:=H^*(S,\mathbb{Z})$ be a lattice, where the pairing is given by intersection on $H^2$ and duality between $H^0$ and $H^4$. We call this lattice the Mukai lattice. It is isometric to $U^4\oplus E_{8}(-1)^2$.
\end{defn}
Markman \cite{mar1} proved that the Mukai lattice can be used to classify Hodge isometries obtained by parallel transport:
\begin{thm}\cite[Theorem 9.3]{mar1}
Let $X$ be a manifold of \kntipo. Then there exists a canonically defined equivalence class of embeddings $\iota_X\colon H^2(X,\mathbb{Z})\rightarrow L_{M}$. A Hodge isometry $g\colon H^2(X,\mathbb{Z})\rightarrow H^2(Y,\mathbb{Z})$ is a parallel transport operator if and only if $\iota_X=\iota_Y\circ g$ in $O(L_n,L_{M})/O(L_{M})$
\end{thm}
For $X=Y$, this amounts to saying that $g$ acts as $\pm1$ on $L_n^\vee/L_n$.

Parallel transport operators can then be used for the following
\begin{thm}[Global Torelli, Huybrechts, Markman and Verbitsky]\label{thm:global_torelli}
Let $X$ and $Y$ be two \hk manifolds of \kntipo. Suppose $\psi\,:\,H^2(X,\mathbb{Z})\,\rightarrow\,H^2(Y,\mathbb{Z})$ is a parallel transport Hodge  isometry. Then there exists a birational map $\phi\,:\,X\,\dashrightarrow\,Y$.
\end{thm}
Related to this there is also the following useful theorem, due to Huybrechts \cite[Theorem 3.2]{mar1}:
\begin{thm}\label{thm:graph}
Let $(X,f)$ be a marked Hyperk\"{a}hler manifold and $(X',g)$ another marked Hyperk\"{a}hler manifold such that the points $(X,f)$ and $(X',g)$ are not separated. 
 Then there exists an effective cycle $\Gamma=Z+\sum_{j}Y_j$ in $X\times X'$ satisfying the following conditions:
\begin{itemize}
\item Z is the graph of a bimeromorphic map from $X$ to $X'$.
\item The codimensions of $\pi_1(Y_j)$ and $\pi_2(Y_j)$ are equal.
\item The composition $g^{-1}\circ f$ is equal to $\Gamma_*\,:\,H^2(X,\mathbb{Z})\,\rightarrow\,H^2(X',\mathbb{Z})$.
\item The cycles $\pi_i(Y_j)$ are uniruled.
\end{itemize}
\end{thm}
Finally, the birational geometry of certain manifolds of \kntiposp has been analysed by Bayer and Macr\'{i} \cite{bm} and was generalised in \cite{me_kahl} (see \cite{bht} for a different approach) using general properties of the \kahl cone of \hk manifolds. 

Let $X$ be a manifold of \kntiposp and let $H^2(X,\mathbb{Z})\rightarrow L_M$ with the primitive embedding described in \cite[Theorem 9.3]{mar1}. Let $v$ be a generator of the orthogonal complement of $H^2$ under this embedding. For every divisor $D$ we denote by  $T_D$ the primitive rank 2 lattice containing $v$ and $D$ inside $L_M$ and we let $v,r_D$ be its generators.
\begin{defn}\label{defn:wall_div}
Let $X,D$ and $T_D$ be as above. The divisor $D$ is a wall divisor if one of the following are satisfied:
\begin{itemize}
\item $r_D^2=-2$ and $0\leq(v,r_D)\leq v^2/2$.
\item $0\leq r_D^2v^2\leq (v,r_D)^2<(v^2/2)^2$.
\end{itemize}

\end{defn}

\begin{thm}\cite[Theorem 1.3 and Proposition 1.5]{me_kahl}
Let $X$ be a manifold of \kntiposp and let $\mathcal{W}$ be the set of wall divisors on $X$. Then the \kahl cone of $X$ is one of the connected components of the following set
\begin{equation}
\{x\in H^2(X,\mathbb{R}),\,x^2>0,\,(x,w)\neq0\,\forall\,w\in\mathcal{W}\}.
\end{equation}
\end{thm}

\section{Lattice theory}\label{sec:lattice}
In this section we sketch some of the lattice theory used in the proofs of the main theorems. An additional list of results, mainly an enumeration of lattices of interest, is contained in the appendix. The interested reader can consult \cite{nik2} for prerequisites on discriminant groups and forms. 
By lattice we mean a free $\mathbb{Z}$ module equipped with a non degenerate bilinear form. We call it even if the associated quadratic form takes only even values. Given an element $v\in L$, we denote $div(v)$ a positive generator of the ideal $(v,L)$ and we call it the divisibility of $v$.
\subsection{Discriminant groups}
For a lattice $L$ its discriminant group is $A_L:=L^\vee/L$. Let $l(A_L)$ denote the length of this group. If the lattice $L$ is even, $A_L$ has a bilinear form with values in $\mathbb{Q}/\mathbb{Z}$ induced from the bilinear form on $L$. Many properties of the associated quadratic form on $\mathbb{Q}/2\mathbb{Z}$, called discriminant form, were found by Nikulin in \cite{nik2}.
If $L$ is a lattice, we call $(l_+,l_-)$ its signature and $q_L$ its discriminant form. The integer $l_+-l_-$ is called signature of $q_L$ and, modulo $8$, it is well defined.\\
Here we will make often use of the following facts concerning primitive embeddings.

\begin{lem}\cite[Proposition 1.15.1]{nik2}\label{lem:nik_immerge}
Let $S$ and $N$ be even lattices. Primitive embeddings of S into $N$ are determined by the sets $(H_S,H_N,\gamma,K,\gamma_K)$, where $K$ is an even lattice with signature  $(n_+-s_+,n_--s_-)$ and discriminant form $-\delta$ where $\delta\,\cong\,(q_{A_S}\oplus -q_{A_N})_{|\Gamma_\gamma^\perp/\Gamma_\gamma}$ 
and $\gamma_K\,:\,q_K\,\rightarrow\,(-\delta)$ is an isometry.\\
Moreover two such sets $(H_S,H_N,\gamma,K,\gamma_K)$ and $(H'_S,H'_N,\gamma',K',\gamma'_K)$ determine isomorphic sublattices if and only if
\begin{itemize}
\item $H_S=\lambda H'_S$, $\lambda\in O(q_S)$,
\item $\exists\,\epsilon\,\in\,O(q_{A_N})$ and $\psi\,\in\,Isom(K,K')$ such that $\gamma'=\epsilon\circ\gamma$ and $\overline{\epsilon}\circ\gamma_K=\gamma'_K\circ\overline{\psi}$, where $\overline{\epsilon}$ and $\overline{\psi}$ are the isometries induced among discriminant groups.
\end{itemize} 
\end{lem}
Here $\Gamma_{\gamma}$ is the graph of $\gamma$.
For most purposes, only the following simplified version will suffice:
\begin{lem}\label{lem:nik_immerge1}
Let $S$ be an even lattice of signature $(s_+,s_-)$. The existence of a primitive embedding of $S$ into some unimodular lattice $L$ of signature $(l_+,l_-)$ is equivalent to the existence of a lattice $M$ of signature $(m_+,m_-)$ and discriminant form $q_{A_M}$ such that the following are satisfied: 
\begin{itemize}
\item $s_++m_+=l_+$ and $s_-+m_-=l_-$.
\item $A_M\cong A_S$ and $q_{A_M}=-q_{A_S}$.
\end{itemize}
\end{lem}
We will also use a result on the existence of lattices, the following is a simplified version of \cite[Theorem 1.10.1]{nik2}
\begin{lem}\label{lem:nik_esiste}
Suppose the following are satisfied:
\begin{itemize}
\item $sign(q_T)\equiv t_+-t_-$ $mod\,8$.
\item $t_+\geq0$, $t_-\geq0$ and $t_++t_-\geq l(A_T)$.
\item There exists a lattice $T'$ of rank $t_++t_-$ and discriminant form $q_T$ over the group $A_T$.
\end{itemize}
Then there exists an even lattice $T$ of signature $(t_+,t_-)$, discriminant group $A_T$ and form $q_{A_T}$.
\end{lem}
The uniqueness of a lattice is usually difficult to prove, but here are two special cases:

\begin{lem}\cite[Corollary 1.13.3]{nik2}\label{lem:nik_unico}
Let $S$ be an even indefinite lattice with signature $(t_+,t_-)$ and discriminant form $q_{A_S}$. Then all lattices with the same signature and discriminant form are isometric to it if $t_++t_-\geq 2+l(A_S)$.

\end{lem}

Let $L$ be a lattice such that $A_L=(\mathbb{Z}_{/(2)})^r$. It is called a $2$-modular lattice and its isometry class is easily determined using its discriminant form. Let us denote by $\Delta$ a discrete invariant, which is $0$ if the discriminant form takes only values $0$ and $1$, and it is $1$ otherwise.
\begin{thm}\cite[Theorem 3.6.2]{nik4}\label{thm:2mod}
A $2$-modular indefinite lattice is uniquely determined by its rank, its signature, the length of its discriminant group and $\Delta$
\end{thm}

\subsection{Isometries, Invariant and Co-invariant Lattices}
In this subsection we analyse two kind of lattices linked to an isometry, namely the co-invariant and invariant lattices. 

\begin{defn}\label{defn:inv_coinv}
Let $R$ be a lattice and let $G\subset O(R)$. Then we define $T_G(R)=R^G$ as the invariant lattice of $G$ and $S_G(R)=T_G(R)^\perp$ as the co-invariant lattice. 
\end{defn}

We will be mainly interested in the following lattices and groups:

\begin{defn}\label{defn:leech_group}
Let $M$ be an even lattice and let $G\subset\,O(M)$. Then $M$ is a Leech type lattice with respect to $G$ if the following are satisfied:
\begin{itemize}
\item $M$ is negative definite.
\item $M$ contains no vectors of square $-2$.
\item $G$ acts trivially on $A_M$.
\item $S_G(M)=M$.
\end{itemize} 
Moreover we call $(M,G)$ a Leech pair and $G$ a Leech type group.
\end{defn}
Leech type lattices of rank at most $11$ have been classified by Nikulin in \cite[Section 1.14]{nik2} and the only such lattice is $E_8(-2)$. The easiest examples of such lattices can be constructed as follows. Let $N$ be a unimodular negative definite lattice with no elements of square $-2$ and let $G\subset O(N)$. Then $(S_G(N),G)$ is a Leech pair.

\begin{lem}\label{lem:G_tors}
Let $R$ be a lattice, and let $G\subset O(R)$. Then the following hold:
\begin{itemize}
\item $T_G(R)$ contains $\sum_{g\in G}gv$ for all $v\in R$.
\item $S_G(R)$ contains $v-gv$ for all $v\in R$ and all $g\in G$.
\item If $R$ is definite then $T_G(R)$ and $S_G(R)$ are nondegenerate.
\item $R/(T_G(R)\oplus S_G(R))$ is of $|G|$-torsion.
\end{itemize}
\begin{proof}
It is obvious that $\sum_{g\in G}gv$ is $G$-invariant for all $v\in R$. For $w\in T_G(R)$ we have $(w,v)=(gw,gv)=(w,gv)$ for all $v\in R$ and all $g\in G$. Therefore $v-gv$ is orthogonal to all $G$-invariant vectors, hence it lies in $S_G(R)$.
Obviously whenever $R$ is definite all of its sublattices are nondegenerate. 
Let $t\in R$, we can write $|G|t=\sum_{g\in G} g(t) + \sum_{g\in G}(t-g(t))$, where the first term lies in $T_G(R)$ and the second in $S_G(R)$.

\end{proof}
\end{lem}

\begin{oss}\label{oss:biggen_g}
Let $L$ be a lattice and let $L\subset L'$ be a primitive embedding into a unimodular lattice. Let $g\in O(L)$ be an isometry which acts trivially on $A_L$. Then there exists an isometry $\overline{g}\,\in\,O(L')$ such that $\overline{g}_{|L}=g$ and $\overline{g}_{|L^\perp}=Id$.
\end{oss}

\begin{defn}\label{defn:k3_lattices}
Let $X$ be a $K3$ surface and let $G\subset Aut(X)$ be a finite abelian group of symplectic morphisms. Nikulin \cite{nik1} proved that the action of $G$ on $H^2(X,\mathbb{Z})$ is uniquely determined. We denote with $S_{G.K3}$ the lattice $S_G(H^2(X,\mathbb{Z}))$. When $G$ is a cyclic group of order $n$, we use $S_{n.K3}:=S_{G.K3}$. For an explicit description of these lattices, the interested reader can consult \cite{hashi} and \cite{gs}.
\end{defn}

\section{From symplectic automorphisms to the Conway group}\label{sec:aut_in_co}

\begin{defn}\label{defn:inv_locus}
Let $X$ be a manifold of \kntiposp and let $G\subset Aut(X)$. We let $T_{G}(X)$ inside $H^2(X,\mathbb{Z})$ be the sublattice fixed by the induced action of $G$ on $H^2(X,\mathbb{Z})$. Moreover we define the co-invariant lattice $S_G(X)\subset H^2(X,\mathbb{Z})$ as $T_G(X)^{\perp}$. 
The fixed locus of $G$ on $X$ will be denoted $X^G$.
\end{defn}

We wish to remark that the map
\begin{equation}\label{numap}
Aut(X)\stackrel{\nu}{\rightarrow} O(H^2(X,\mathbb{Z}))
\end{equation}
is injective for manifolds of \kntipo.
We have the following exact sequence for any finite group $G$ of Hodge isometries on $H^2(X,\mathbb{Z})$:
\begin{equation}\label{exactgroup}
1\,\rightarrow\, G_0\,\rightarrow\,G\,\stackrel{\pi}{\rightarrow}\,\Gamma_m\,\rightarrow\,1,
\end{equation}
where $\Gamma_m\subset U(1)$ is a cyclic group of order $m$. In fact the action of $G$ on $H^{2,0}$ is the action of a finite group on $\mathbb{C}$.

\begin{lem}\label{lem:gaction_gen}
Let $X$ be a manifold of \kntiposp and let $G\subset Aut(X)$ be a finite group. Then the following hold:
\begin{enumerate}
\item $g\in G$ acts trivially on $T(X)$ if and only if $g\in G_0$.
\item The representation of $\Gamma_m$ on $T(X)\otimes\mathbb{Q}$ splits as the direct sum of irreducible representations of the cyclic group $\Gamma_m$ having maximal rank (\ie of rank $\phi(m)$).
\end{enumerate}
\begin{proof}
The proof goes exactly as \cite[Lemma 3.4]{me_invol}. See \cite[Proposition 6]{beau2} or \cite[Lemma 7.1.4]{me_tesi} for further reference.
\end{proof}
\end{lem}

\begin{defn}
Let $G\subset O(L_n)$ be a group and let $X$ be a manifold of \kntiposp such that $Pic(X)\cong S_G(L_n)$. Suppose moreover that $G$ consists of parallel transport operators on $X$. A numerical wall divisor in $S_G(L_n)$ is the image under the above isometry of a wall divisor of $X$.
\end{defn}

\begin{lem}\label{lem:discr_preserving}
Let $f\in O(L_n)$ be an isometry such that $f$ acts as $-Id$ on $A_{L_n}$ and $S_f(L_n)$ is negative definite. Then $S_f(L_n)$ contains a numerical wall divisor.
\begin{proof}
Isometries of $L_n$ are an extension of discriminant preserving isometries with isometries of $A_{L_n}$, therefore we can write $f=g\circ R_v$, where $R_v$ is the reflection along an element of square $2-2n$ and $div(v)=2n-2$. In particular we have that $S_g(L_n)$ is also negative definite. Let $f(v)=-v+(2n-2)w$. This implies $v-g(v)=(2-2n)w$, \ie $w\in S_g(L_n)$. Moreover $t:=(v-f(v))/2=v-(n-1)w\in S_f(L_n)$ is a non zero element. Since $f$ is an isometry, $(v,w)=(n-1)w^2\leq 0$ because $w\in S_g(L_n)$. Therefore $2-2n\leq t^2<0$ and $div(t)\geq n-1$. Let us now choose any manifold $X$ with $Pic(X)\cong S_f(L_n)$ such that $f$ preserves the embedding $H^2(X,\mathbb{Z})\rightarrow L_M$. Let $s$ be a generator of $(H^2)^\perp$ in this embedding. The reflection by $v$ is a parallel transport operator, hence either $\frac{s+t}{div(t)}$ or $\frac{s-t}{div(t)}$ lie in $L_M$. In any case, the primitive lattice containing $s$ and $t$ satisfies \Ref{defn}{wall_div}, therefore $t$ is a numerical wall divisor.
\end{proof}
\end{lem}
We remark that the above condition on the signature of $S_{f}(L_n)$ is indeed necessary, there are examples of nonsymplectic automorphisms acting nontrivially on the discriminant group of the manifold. 
Let now $G$ be a finite group of automorphisms of $X$ such that $G\subset Aut_s(X)$.
\begin{lem}\label{lem:algaction}
Let $X$ be a \hk manifold and let $G\subset Aut_s(X)$ be a finite group. Then the following assertions are true:
\begin{enumerate}
\item $S_G(X)$ is nondegenerate and negative definite.
\item $T(X)\subset T_G(X)$ and $S_G(X)\subset S(X)$.
\item $S_G(X)$ contains no wall divisors
\item The action of $G$ on $A_{S_G(X)}$ is trivial.
\end{enumerate}
\begin{proof}
The proof of the first three items goes exactly as in \cite[Lemma 3.5]{me_invol} with wall divisors taking the role of -2 elements. For the reader's convenience we sketch it here. The invariant lattice $T_G(X)$ contains $T(X)$ because $G$ is symplectic and, after tensoring with $\mathbb{R}$, it contains an invariant \kahl class because $G$ is finite. Therefore its orthogonal $S_G(X)$ is negative definite. Since $T_G(X)\otimes \mathbb{R}$ contains a \kahl class, its orthogonal can not contain wall divisors.\\  
For the final statement, suppose that an element $g$ of $G$ acts nontrivially on $A_{L_n}$. Since it is a parallel transport operator, it must act as $-Id$ and \Ref{lem}{discr_preserving} implies that $S_g(X)$ would contain a wall divisor.
\end{proof}
\end{lem}

We are now ready to prove the main result of this section:
\begin{thm}\label{thm:sporadic}
Let $X$ be a manifold of \kntiposp and let $G\subset Aut_s(X)$ be a finite group. Then $G\subset Co_1$ and $S_G(X)=S_G(\Lambda)$ for some conjugacy class of $G$ in $Co_1$. 
\begin{proof}
The first part of the theorem is a generalisation of \cite[Theorem 1.1]{me_autom1}, however we provide a different proof which simplifies the classification and yields also the second part of the theorem. This proof is essentially taken from \cite{huy_der}.\\ 
Let us take a primitive embedding $S_G(X)\rightarrow \Lambda\oplus U$, which exists by \Ref{lem}{nik_immerge1}. By \Ref{lem}{discr_preserving}, the action of $G$ can be extended by $Id$ to an action on $\Lambda\oplus U$. Let us denote with $e$ one of the standard generators of the above copy of $U$. Let $V=(\Lambda\oplus U)\otimes\mathbb{R}$ and let us consider a decomposition of its positive cone in a wall and chamber structure, where walls are orthogonal to elements of square $-2$. Then $T_G(V)$ intersects one of these chambers, otherwise we would have an element of square $-2$ contained in its orthogonal. The action of the Weyl group of reflections is transitive on these chambers, therefore we may assume that $G$ preserves a special chamber, denoted $\mathcal{C}_0$, which is cut out by all roots $v$ such that $(v,e)=1$. The group of isometries fixing $\mathcal{C}_0$ is $Co_\infty$, see \cite[Chapter 27]{con}. The group $G$ also fixes $e$, therefore we have an embedding $S_G(X)\rightarrow \Lambda$. The action of $G$ on $A_{S_G(X)}$ is trivial, therefore $G$ extends to a group of isometries of $\Lambda$. Furthermore we have $G\subset Aut(\Lambda)/\{\pm1\}=Co_1$ because $-Id$ has a coinvariant lattice of rank $24$. 
\end{proof}
\end{thm}

\section{From the Conway group to symplectic automorphisms}\label{sec:co_in_aut}

\begin{thm}\label{thm:cohom_to_aut}
Let $G$ be a finite subgroup of $O(L_n)$. Then $G$ is induced by a symplectic subgroup of $Aut(X)$ for some manifold $(X,f)$ of \kntiposp if and only if the following hold:
\begin{enumerate}
\item $S_G(L_n)$ is non degenerate and negative definite.
\item $S_G(L_n)$ contains no numerical wall divisor.
\end{enumerate}

\begin{proof}
By the surjectivity of the period map and by \Ref{lem}{algaction} we can consider a marked manifold $(X,f)$ of \kntiposp such that $T(X)\stackrel{f}{\rightarrow} T_G(L_n)$ is an isomorphism and also $S(X)\stackrel{f}{\rightarrow} S_G(L_n)$ is.\\
Since $S_G(L_n)$ contains no numerical wall divisors, it follows that $Pic(X)$ contains no wall divisors and therefore $\mathcal{K}_X=\mathcal{BK}_X=\mathcal{C}_X$. In particular we have a $G$ invariant \kahl class and $X$ contains no effective rational curves nor divisors.\\
For $g\in G$, we consider the marked varieties $(X,f)$ and $(X,g\circ f)$. Since $g$ acts trivially on the discriminant group of $X$, these two marked manifolds lie in the same connected component of the moduli space. They also have the same period, hence by \Ref{thm}{graph} we have $f^{-1}\circ g\circ f=\Gamma_*$. Here $\Gamma=Z+\sum_{j}Y_j$ in $X\times X$, where $Z$ is the graph of a bimeromorphic map from $X$ to itself and $Y_j$'s are cycles with $codim(\pi_i(Y_j))\geq 1$.\\ 
 All $Y_j$'s contained in $\Gamma$ have $codim(\pi_i(Y_j))>1$, thus implying $g=\Gamma_*=Z_*$ on $H^2(X,\mathbb{Z})$. Now the bimeromorphic map $Z$ is biregular since $\mathcal{K}_X=\mathcal{C}_X$. 
\end{proof}
\end{thm}

\begin{oss}
A wall divisor $D$ on a manifold $X$ of \kntiposp has either square $-2$ or $div(D)>1$ and it divides $2(n-1)$. So, if $G$ is a cyclic group of order coprime with $2(n-1)$, the only numerical wall divisors that $S_G(L_n)$ could contain are those of square $-2$.
\end{oss}

\begin{prop}\label{prop:sfiga}
Let $(S,G)$ be a pair consisting in a Leech type lattice and its Leech automorphism group as in \Ref{defn}{leech_group}. Let moreover $S\subset N$, one of the 24 Niemeier lattices.\\ Suppose there exists a primitive embedding $S\rightarrow L_n$ and suppose that the order of $G$ is coprime with $2(n-1)$.\\
Then $G$ extends to a group of automorphisms on some manifold $X$ of \kntipo.
\begin{proof}
This is an immediate consequence of \Ref{thm}{cohom_to_aut}: $G$ acts trivially on $A_S$, therefore we can extend $G$ to a group of isometries of $L_n$ acting trivially on $S^{\perp_{L_n}}$. Thus we have $S_G(L_n)\cong S$. Moreover since $S$ is a Leech type lattice contained in a negative definite lattice $N$ the other conditions of \Ref{thm}{cohom_to_aut} are satisfied.
\end{proof}
\end{prop}

\begin{thm}\label{thm:da_conway_ad_auto}
Let $G\subset Co_0$ be a group of isometries such that $rk(S_G(\Lambda))\leq 20$ and $rk(T_G(\Lambda))>l(A_{T_G(\Lambda)})$. Then there exist an integer $n$ and a manifold $X$ of \kntiposp such that $G\subset Aut_s(X)$ and $S_G(X)\cong S_G(\Lambda)$.
\begin{proof}
First of all, by the existence of the lattice $T_G(\Lambda)$, we have a lattice $T$ of signature $(4,rk(T_G(\Lambda))-4)$ and $q_T=q_{A_{T_G(\Lambda)}}$ by \Ref{lem}{nik_esiste}. Then, by \Ref{lem}{nik_immerge1}, there is a primitive embedding of $S_G(\Lambda)$ into the Mukai lattice $L_M$ whose orthogonal is $T$. We have $rk(T)>l(A_T)$. This implies that there exists an element $v\in T$ such that $(v+t)/r\in L_M$ implies $r=1$ for all $t\in S_G(\Lambda)$. After adding with some positive element of $T$, we have $v^2\geq2$ (and it can actually take infinitely many positive values). Let now $n=(v^2+2)/2$, the lattice $v^\perp$ is isometric to $L_n$ and all elements of $S_{G}(\Lambda)$ have divisibility $1$ in $L_n$. Since $S_G(\Lambda)$ contains no elements of square $-2$ and no elements of divisibility at least $2$, it contains also no numerical wall divisors. Therefore \Ref{thm}{cohom_to_aut} applies and we obtain our claim.
\end{proof}
\end{thm}
The condition in the above theorem is actually equivalent to $rk(S_G(\Lambda))+l(A_{S_G(\Lambda)})\leq 23$. This equivalence is trivial if $rk(S_G(\Lambda))\geq 12$ and it still holds otherwise because the only smaller Leech type lattice is $E_8(-2)$.

With the above proposition, we can construct automorphism groups from several subgroups of the Conway group. In \cite[Theorem 0.2]{huy_der}, Huybrechts also finds a sufficient condition to obtain a group of symplectic automorphisms on a manifold of \kntiposp from a subgroup of the Conway group, the following shows that his condition is equivalent to ours:
\begin{lem}\label{lem:cond_huy}
Let $M$ be a negative definite lattice of rank at most $20$. Then the following are equivalent:
\begin{itemize}
\item $M$ embeds primitively in $\Lambda$ and its orthogonal $T_{\Lambda}$ satisfies $rk(T_{\Lambda})> l(A_{T_\Lambda})$. 
\item $M$ embeds primitively in $L_M$ and its orthogonal $R$ satisfies $rk(R)> l(A_{R})$.
\item $M$ embeds primitively into a lattice $P$ of signature $(1,20)$ with $l(A_P)\leq 2$.
\item $M$ embeds primitively in $L_M$ and its orthogonal $R$ contains a positive definite lattice $\Gamma$ of rank $3$ with $l(\Gamma)\leq 2$.
\end{itemize}
\begin{proof}
The last two conditions are trivially equivalent, after noticing that both $\Gamma$ and $P$ embed into $L_M$ and taking $\Gamma=P^\perp$.\\
The first two conditions are also equivalent by \Ref{lem}{nik_esiste} and \Ref{lem}{nik_immerge1}. Let now $\Gamma\subset R$ with $l(\Gamma)\leq 2$. This means that there exists an element $t\in \Gamma$ such that $t/p\notin \Gamma^{\vee}$ for any prime $p$ and also $t/p\notin R^\vee$, \ie $l(A_R)< rk(R)$. Equivalently if $l(A_R)< rk(R)$ there exists a $t$ as above and we define $\Gamma$ as any lattice containing $t,v $ and $w$ where these three are linearly independent.
\end{proof}
\end{lem}

One is naturally interested in the case of a group $G\subset Aut(\Lambda)$ such that $rk(T_G(\Lambda))=l(A_{T_G(\Lambda)})$. The following proposition actually shows that, if ever we could obtain $G$ as an automorphism of manifolds of \kntipo, this would happen only for finitely many $n$. 

\begin{prop}\label{prop:wall_in_s}
Let $G\subset O(L_n)$ be a group such that $M:=S_G(L_n)<0$ and $G$ acts trivially on $A_{L_n}$. Suppose moreover that $l(A_M)+rk(M)=24$. Let $\{t\}_i$ be a set of primitive elements of $M$ such that $[t_i/div(t_i)]$ are the nontrivial part of $A_M$ and $|t_i^2|\leq div(t_i)^2(n+3)/2$. Then $M$ contains a wall divisor for any manifold $X$ with $M\subset Pic(X)$.
\begin{proof}
Let $L_n\subset L_M$ and let $v$ be a generator of $L_n^\perp$. Let $T:=M^\perp$ in $L_M$. The condition $l(A_M)+rk(M)=24$ is equivalent to $rk(T)=l(A_T)$. Therefore, there exists an element $t\in M$ such that $(v+t)/a\in L_M$, where $a=GCD(div_{M}(t),2n-2)\neq 1$. Moreover, $t=t_i+div(t)w$ for some $i$ and some $w\in M$. Let now $X$ be a manifold of \kntiposp with $M\subset Pic(X)$ and the embedding $L_n\subset L_M$ coincides with the natural one on $X$. By our assumption on $t_i$, the vector $(v+t_i)/a$ is a wall divisor as in \Ref{defn}{wall_div}, hence we are done.  
\end{proof}

\end{prop}
In particular, if we analyse what happens for three interesting lattices with $l(A_T)+rk(T)=24$, we see that they never occur as coinvariant lattices for symplectic automorphisms. This is some evidence for the following:
\begin{conj}
There is a bijective correspondence between finite groups of symplectic automorphisms of manifolds of \kntiposp (for some $n$) and subgroups $G$ of $Co_1$ satisfying the conditions of \Ref{thm}{da_conway_ad_auto}.
\end{conj} 

The three lattices we will analyse are coinvariant lattices for isometries of order $2$ or $3$ and they are the Barnes-Wall lattice $BW_{16}(-1)$, defined in \cite[Section 4.10]{con}, the lattice $S_{3.exo}$ as in \Ref{ex}{p3E83} and the lattice $D^+_{12}(-2)$, where $D_{12}^+$ is an (odd) unimodular overlattice of $D_{12}$.
\begin{lem}\label{lem:bw16_embed}
Let $v\in L_M$ be a primitive element of square $2n-2$ and let $BW_{16}(-1)$ be primitively embedded inside $L_n:=v^\perp$. Then there exists an element $t$ inside $BW_{16}$ such that $\frac{v+t}{2}\in L_M$ and $t$ is a numerical wall divisor.
\begin{proof}
We want to apply \Ref{prop}{wall_in_s}, and we use the fact that the discriminant group of $BW_{16}$ can be generated by elements $t_i$ of square at most $12$, as computed in \cite[Section 6.5]{con}. Moreover, by \Ref{thm}{2mod}, the orthogonal of $BW_{16}$ inside $L_M$ is isometric to $U(2)^4$ hence $n$ must be odd and the inequality of \Ref{prop}{wall_in_s} is satisfied.
\end{proof}
\end{lem}
\begin{lem}\label{lem:s3exo_embed}
Let $v\in L_M$ be a primitive element of square $2n-2$, where $n\equiv 1$ mod $3$. Let $S_{3.exo}$ be primitively embedded inside $L_n:=v^\perp$. Then there exists an element $t$ inside $S_{3.exo}$ such that $\frac{v+t}{3}\in L_M$ and $t$ is a numerical wall divisor.
\begin{proof}
The proof goes as in \Ref{lem}{bw16_embed}. Here the set $\{t_i\}$ is given by elements of the form $2e-f-g$, where $S_{3.exo}\subset E_8(-1)^3$, $e$ is a root of one copy of $E_{8}$ and $f$ and $g$ are the corresponding roots in the other two copies. These have squares $-36,-24$ or $-12$ which satisfies the bound in \Ref{prop}{wall_in_s} if $n\geq 5$. Observing that $t^2$ and $2n-2$ must be congruent modulo $9$, we are done also for $n\leq 4$.

\end{proof}
\end{lem}

\begin{lem}\label{lem:d12_embed}
Let $v\in L_M$ be a primitive element of square $2n-2$ and let $D^+_{12}(-2)$ be primitively embedded inside $L_n:=v^\perp$. Then there exists an element $t$ inside $D^+_{12}$ such that $\frac{v+t}{2}\in L_M$ and $\langle v,\frac{v+t}{2}\rangle$ is primitive and isometric to $\left(\begin{array}{cc} 2-2n & n-1\\ n-1 & -2\end{array}\right)$.
\begin{proof}

We will apply \Ref{lem}{nik_immerge} to the lattice $D^+_{12}(-2)$ several times and we use the fact that there exists an isometry between isometric elements of $A_{D^+_{12}(-2)}$. 
First of all, up to isometry there is only one primitive embedding of $D^+_{12}(-2)$ inside $L_M$, whose orthogonal is $T:=U(2)^4\oplus(-2)^4$ by the classification of $2$-modular lattices in \Ref{thm}{2mod}. Let $v\in T$ be as above, the lattice spanned by $v$ and $D^+_{12}(-2)$ is primitive inside $L_M$ and is isometric to $U\oplus (-2)^{11}$ again by \Ref{thm}{2mod}.\\ Once more, also the primitive embedding of $D^+_{12}(-2)$ inside $U\oplus (-2)^{11}$ is unique up to isometry. One such embedding is given by taking an element $t$ of square $-2n-6$ and adding $\frac{v+t}{2}$ to $v\oplus D^+_{12}(-2)$, which is our claim. 
\end{proof}
\end{lem}

\section{Classification of prime order morphisms}\label{sec:prime}

The aim of this section is to give a first application of \Ref{thm}{sporadic}, \ie the classification of prime order symplectic automorphisms on manifolds of \kntiposp up to their co-invariant lattice and the minimal $n$ for which they occur, the case $n=1$ being morphisms of $K3$ surfaces. 

\begin{cor}\label{cor:k3n_prime}
Let $X$ be a manifold of \kntiposp and let $\varphi\in Aut_s(X)$ be of prime order $p$. Then one of the following holds:
\begin{table}[ht]\label{tab:prime_autom_k3n_tab}
\begin{tabular}{|c|c|c|}
\hline
$p$ & Lattice $S_{\varphi}(X)$ & Minimal $n$\\
\hline
2 & $E_8(-2)$ & 1\\
\hline
3 &  $S_{3.K3}$, as in \Ref{defn}{k3_lattices} & $1$\\
\hline
3 &  $W(-1)$ as in \Ref{ex}{slat_rk6} & $2$\\
\hline
5 &  $S_{5.K3}$ as in \Ref{defn}{k3_lattices} & $1$\\
\hline
5 &  $S_{5.exo}$ as in \Ref{ex}{slat_5rk4} & $3$\\
\hline
7 &  $S_{7.K3}$ as in \Ref{defn}{k3_lattices} & $1$\\
\hline
11 & $S_{11.K3^{[2]}}$ as in \Ref{ex}{p11A212} & $2$\\
\hline
\end{tabular}
\end{table}
\begin{proof}
The list of lattices is just a direct consequence of \Ref{prop}{invol_nieme}, \Ref{prop}{p3_nieme}, \Ref{prop}{p5_nieme}, \Ref{prop}{p7_nieme}, \Ref{ex}{p11A212} and \Ref{prop}{no13}. All these lattices satisfy the conditions of \Ref{prop}{sfiga}, therefore they are covariant lattices of symplectic automorphisms. 
We have only to exclude the cases $S_\varphi=D^+_{12}(-2),\,BW_{16}$ or $S_{3.exo}$. All these lattices embed in the Mukai lattice, as shown in \Ref{lem}{d12_embed}, \Ref{lem}{bw16_embed} and in \Ref{lem}{s3exo_embed}. However, all these embeddings contain a numerical wall divisor. Therefore the Leech isometry of these lattices does not induce an automorphism.\\ Obviously all the cases corresponding to natural automorphisms exist in all possible dimensions.
To analyse all the other cases we must embed the lattices $S_i$ contained in \Ref{tab}{prime_autom_k3n_tab} inside the Mukai lattice $L_M$ and look at their orthogonal: if it represents the integer $2(n-1)$ with a primitive vector then there exists a primitive embedding  of $S_i$ inside $L_n$. A little care is required if $p$ divides $n-1$, since in that case some embeddings of $S_i$ might contain numerical wall divisors. Otherwise by \Ref{thm}{cohom_to_aut} there exists a manifold of \kntiposp having an automorphism $\varphi$ such that $S_{\varphi}\cong S_i$. Let us look at all cases one by one:
\begin{itemize}
\item[$S_i=W(-1)$] In this case we will look at a greater lattice: let $F$ be the orthogonal inside $\Lambda$ to the $\mathcal{S}$-lattice (see \Ref{defn}{slat}) $2^93^6$. Then $W(-1)\subset F$ by $2^93^6\subset 2^{27}3^{36}$. Let us now embed $F$ into $L_M$ and let $T=F^{\perp_L}$. A necessary condition for $W(-1)\rightarrow L_n$ is that $2(n-1)$ is represented by a primitive vector of $T$. By \Ref{ex}{slat_3rk4}, $T\cong A_2\oplus A_2(3)$. Therefore it represents $2$, so $W(-1)$ primitively embeds into $L_2$.
\item[$S_i=S_{5.exo}$] Let us fix an embedding $S_{5.exo}\rightarrow L_M$ and let $T=S_{5.exo}^\perp$. By \Ref{ex}{slat_5rk4}, $T$ is in the same genus of the $\mathcal{S}$-lattice $2^53^{10}(-1)$, however there is only one lattice in this genus, which we recall is
\begin{equation*}
\left(
\begin{array}{cccc}
4 & 1 & 1 & -1\\
1 & 4 & -1 & 1\\
1 & -1 & 4 & 1\\
-1 & 1 & 1 & 4
\end{array}
\right).
\end{equation*}
An easy computation shows that the minimal integer it represents is $4$,  therefore $S_i$ does embed into $L_3$. 
\item[$S_i=S_{11.K3^{[2]}}$] Let us fix an embedding $S_{11.K3^{[2]}}\rightarrow L_M$ and let $T=S_{11.K3^{[2]}}^\perp$. $T$ has determinant 121 and, as shown in \cite{nip}, there is only one genus of such lattices, containing the following:
\begin{equation*}
\left(
\begin{array}{cccc}
4 & 2 & 1 & 0\\
2 & 4 & 1 & 1\\
1 & 1 & 4 & 2\\
0 & 1 & 2 & 4
\end{array}
\right),\,\,\,
\left(
\begin{array}{cccc}
2 & 1 & 1 & 0\\
1 & 2 & 1 & 1\\
1 & 1 & 8 & 4\\
0 & 1 & 4 & 8
\end{array}
\right),\,\,\,
\left(
\begin{array}{cccc}
2 & 0 & 1 & 0\\
0 & 2 & 0 & 1\\
1 & 0 & 6 & 0\\
0 & 1 & 0 & 6
\end{array}
\right).
\end{equation*}
A direct computation shows that the integer $2$ is represented by two of these lattices, therefore $S_{11.K3^{[2]}}$ embeds into $L_2$ in two different ways.
\end{itemize}

\end{proof}
\end{cor}
In the case of dimension 4, we can be more precise and determine also the fixed locus $X^\varphi$ and the number of deformation types of pairs $(X,\varphi)$.
\begin{cor}\label{cor:k32_prime}
Let $X$ be a manifold of \ktiposp and let $\varphi\in Aut(X)$ be a symplectic morphism of prime order $p$. Then one of the following occurs
\begin{table}[ht]\label{tab:prime_autom_tab}
\begin{tabular}{|c|c|c|c|}
\hline
$p$ & Fixed locus $X^{\varphi}$ & Lattice $S_{\varphi}(X)$ & Deformation classes\\
\hline
2 & $1$ K3 surface and 28 isolated points & $E_8(-2)$ & 1\\
\hline
3 & 27 isolated points & $S_{3.K3}$ & 1\\
\hline
3 & 1 abelian surface & $W(-1)$  & 1\\
\hline
5 & 14 isolated points & $S_{5.K3}$  & 1\\
\hline
7 & 9 isolated points & $S_{7.K3}$ & 1\\
\hline
11 & 5 isolated points & $S_{11.K3^{[2]}}$  & 2\\
\hline
\end{tabular}
\end{table}
\begin{proof}
The list of lattices is given just by the elements of \Ref{cor}{k3n_prime} occurring for $n=2$. To determine the deformation classes of pair, we can use the main theorem of \cite{me_natural}. Indeed, an identical technique can be used to deform two pairs $(X,\varphi)$ and $(Y,\psi)$ if they have isometric invariant and coinvariant lattices. For all cases arising as morphisms of a $K3$ surface, we have uniqueness by \Ref{lem}{nik_unico}. In the case of order $11$, there are two choices for the lattice $T_{\varphi}(X)$ and therefore two deformation classes of pairs. These have been made explicit in \cite[Examples 4.5.1 and 4.5.2]{me_tesi}. Finally, the deformations with coinvariant lattice $W(-1)$ can be analysed as in \Ref{cor}{k3n_prime} by taking a bigger group of automorphisms containing $\varphi$, namely the group given as an extension of $(\mathbb{Z}/3\mathbb{Z})^4$ by $A_6$. Here we have only one deformation class and we can always deform a manifold with $Pic\cong W(-1)$ to one with $Pic\cong F$ by adding classes without adding wall divisors. Once we have an example in each deformation class, we just need to compute the fixed locus in a single example since it deforms smoothly. Such examples are classically known and explicit computations can be found in \cite[Chapter 4]{me_tesi}.
\end{proof}
\end{cor}

\appendix
\section{Sublattices of the Leech lattice}  
\subsection{$\mathcal{S}$-lattices}
In this subsection we analyse briefly a few sublattices of the Leech lattice which arise as $T_G(\Lambda)$ for some interesting groups $G$. Let us start with the basics:
\begin{defn}\label{defn:slat}
Let $M\subset \Lambda$. Then $M$ is a $\mathcal{S}$-lattice \index{Lattice, Special sublattice of $\Lambda$, $\mathcal{S}$-lattice} if all elements of $M$ are congruent modulo $2\Lambda$ to an element of $M$ of norm $0,-4$ or $-6$.
\end{defn}
There are not many examples of $\mathcal{S}$-lattices and they were classified by Curtis:
\begin{lem}\cite{cur}
Up to isomorphisms there are 12 $\mathcal{S}$-lattices inside $\Lambda$.
\end{lem}
Their stabilizers and their automorphism groups inside $Co_0$ are also classified, a full table can be found in \cite[page 180]{atlas}.
For our purpose it is better to give an explicit presentation of the Leech Lattice $\Lambda$:
\begin{ex}\label{ex:leech_24}
Let us consider the vector space $\mathbb{R}^{W}$, where $W=\mathbb{P}^1(\mathbb{Z}_{/(23)})$ is a set with 24 elements and let us endow it with  a quadratic form defined as the opposite of the euclidean norm. Let $Q\subset W$ be the set whose elements are quadratic residues modulo $23$ and $0$, and let $a=8^{-1/2}$. Then $\Lambda\subset \mathbb{R}^W$ is spanned by the following elements:
\begin{align*}
&a(2,\dots,2,0,\dots,0), & \text{    where the twelve non zero elements are}\\
&& \text{supported on a translate of $Q$ by an element of $W$},\\
&a(-3,1,\dots,1), & \\
&a(\pm4,\pm4,0,\dots,0). & 
\end{align*}  
\end{ex} 
Let us introduce a piece of notation: a $\mathcal{S}$-lattice $M$ is called $2^i3^j$ if (up to sign) it contains $i$ vectors of norm $-4$ and $j$ vectors of norm $-6$.
\begin{ex}
The easiest example possible is that of a lattice $M=(-4)=2^1$ in the above notation. The condition of \Ref{defn}{slat} is trivially satisfied and $Aut(M)=\pm Id$, $Stab(M)=Co_2$ \ie $M=T_{Co_2}(\Lambda)$.
\end{ex} 

\begin{ex}\label{ex:slat_5rk4}
Let us consider the $\mathcal{S}$-lattice $M=2^53^{10}$, it has rank 4 and it is $T_G(\Lambda)$, where $G$ is an extension of $(\mathbb{Z}_{/(5)})^3$ with $\mathbb{Z}_{/(4)}$. We wish to remark that $G$ contains an element of the conjugacy class $5C$ in the notation of \cite{atlas}. We will later denote $M^\perp$ as the lattice $S_{5.exo}$ \index{Lattice, $S_{5.exo}$}, which, as shown in \Ref{ex}{A46}, is $S_{5C}(\Lambda)$. Moreover $M$ is isometric to the following
\begin{equation*}
\left(
\begin{array}{cccc}
-4 & -1 & -1 & 1\\
-1 & -4 & 1 & -1\\
-1 & 1 & -4 & -1\\
1 & -1 & -1 & -4
\end{array}
\right).
\end{equation*} 
From Nipp's \cite{nip} list of definite quadratic forms we have that $M$ is the unique lattice in its genus, moreover $M\oplus M$ has $E_8(-1)$ as an overlattice. 
\end{ex}

\begin{ex}\label{ex:slat_3rk4}
The $\mathcal{S}$-lattice $M=2^93^6$ is a lattice of rank 4 and it is the stabilizer of a group $G\subset Co_0$ which is a nontrivial extension of $(\mathbb{Z}_{/(3)})^4$ with $A_6$. If we consider $\Lambda$ as the lattice defined in \Ref{ex}{leech_24} then $M$ is spanned by the following 9 elements:
\begin{align*}
&a(0,0,0,0,0,0,0,0,-4,0,0,0,0,0,0,0,4,0,0,0,0,0,0,0), & \\
&a(4,0,0,0,0,0,0,0,0,0,0,0,0,0,0,0,-4,0,0,0,0,0,0,0), & \\
&a(-4,0,0,0,0,0,0,0,4,0,0,0,0,0,0,0,0,0,0,0,0,0,0,0), & \\
&a(0,0,0,0,0,0,0,0,2,2,2,2,0,0,0,0,-2,-2,-2,-2,0,0,0,0), & \\
&a(-2,-2,-2,-2,0,0,0,0,0,0,0,0,0,0,0,0,2,2,2,2,0,0,0,0), &\\
&a(2,2,2,2,0,0,0,0,-2,-2,-2,-2,0,0,0,0,0,0,0,0,0,0,0,0), &\\
&a(0,0,0,0,0,0,0,0,2,-2,-2,-2,0,0,0,0,-2,2,2,2,0,0,0,0), & \\
&a(-2,2,2,2,0,0,0,0,0,0,0,0,0,0,0,0,2,-2,-2,-2,0,0,0,0), & \\
&a(2,-2,-2,-2,0,0,0,0,-2,2,2,2,0,0,0,0,0,0,0,0,0,0,0,0), &
\end{align*}
where $a=8^{-1/2}$. A direct computation shows that it is isometric to the lattice 
\begin{equation*}
\left(\begin{array}{cccc} 
-4 & 2 & -2 & 1 \\
2 & -4 & 1 & -2 \\
-2 & 1 & -4 & 2 \\
1 & -2 & 2 & -4 \\
\end{array}\right).
\end{equation*}
A look at Nipp's table \cite{nip} shows again that it is unique in its genus. This time however $M\oplus M$ has not an unimodular overlattice, however $M\oplus A_2(-1)\oplus A_2(-3)$ does. Notice moreover that $A_2(-1)\oplus A_2(-3)$ is again unique in its genus by \cite{nip}.
\end{ex}

\begin{ex}\label{ex:slat_rk6}
The $\mathcal{S}$-lattice $M=2^{27}3^{36}$ is a lattice of rank 6 and discriminant group $(\mathbb{Z}_{/(3)})^5$ and it is the stabilizer of a group $G\subset Co_0$, where $G$ is a nontrivial extension of $(\mathbb{Z}_{/(3)})^5$ with $\mathbb{Z}_{/(2)}$. Its orthogonal inside $\Lambda$ is a lattice which contains the group $O(E_6)$ in its automorphism group. In particular, its orthogonal is the Wall lattice $W(-1)$ described in \cite{wal}. This Wall lattice is the coinvariant lattice of an order three isometry of the Leech lattice.
\end{ex}

\subsection{Niemeier lattices}\label{sec:niemeier}

Now we recall Niemeier's list of definite even unimodular lattices of dimension 24. Usually they are defined as positive definite lattices. For our purposes, we will consider them as negative definite lattices. All of these lattices can be obtained by specifying a 0 or 24 dimensional Dynkin diagram such that every semisimple component has a fixed Coxeter number, in \Ref{tab}{nieme} we recall the possible choices. Having the Dynkin lattice $A(-1)$ of the lattice $N$, we obtain $N$ by adding a certain set of glue vectors, which are a subset $G(N)$ of $A^{\vee}/A$. The precise definition of the glue vectors can be found in \cite[Chapter 16]{con} and we keep the same notation contained therein. Notice that the set of glue vectors forms an additive subgroup of $A^{\vee}/A$.\\
Another fundamental data is what we call maximal Leech-type group $Leech(N)$, \ie the maximal subgroup $G$ of $Aut(N)$ such that $(S_G(N),G)$ is a Leech-type pair. It is a well known fact that this group is obtained as $Aut(N)/W(N)$, where $W(N)$ is the Weyl group generated by reflections on $-2$ vectors. These groups where first computed by Erokhin \cite{ero}.\\
This data is summarized in \Ref{tab}{nieme}, taken from \cite{me_autom1}. Let us explain briefly the notation used therein: for the Leech-type group we used standard notation from \cite{atlas}, where $n$ denotes a cyclic group of order $n$, $p^n$ \index{Group, Elementary $p$-group of order $p^n$, $p^n$} denotes an elementary $p$-group of order $p^n$, $G.H$ denotes any group $F$ with a normal subgroup $G$ such that $F/G=H$ and $L_m(n)$ denotes the group $PSL_m$ over the finite field with $n$ elements. $M_n$ denotes the Mathieu group on $n$ elements and $Co_n$ denotes Conway groups.\\   
Regarding the glue codes we kept the notation of \cite{con}, hence a glue code $[abc]$ means a vector $(g,h,f)$ where $g$ is the glue vector of type $a$, $h$ is the one of type $b$ and $f$ of type $c$. Moreover $[(abc)]$ indicates all glue vectors obtained from cyclic permutations of $\{a,b,c\}$, hence $[abc],[bca],[cab]$.\\
\begin{table}[ht]\label{tab:nieme}
\caption{Niemeier lattices and their Leech automorphisms}
\begin{tabular}{|c|c|c|c|c|}
 \hline
Name & Dynkin & Leech-type& Coxeter& Generators of the glue code\\ 
 & diagram & Group &  Number & \\

\hline
$N_{1}$\index{Lattice, Niemeier, $N_i$} & $D_{24}$ & $1$ & $46$ & $[1]$\\ \hline
$N_{2}$ &$D_{16}E_8$ & $1$ & $30$ & $[10]$\\ \hline
$N_{3}$ &$E_8^3$ & $S_3$ & $30$ & $[000]$\\ \hline
$N_{4}$ &$A_{24}$ & $2$ & $25$ & $[5]$\\ \hline
$N_{5}$ &$D_{12}^2$ & $2$ & $22$ & $[12],[21]$\\ \hline
$N_{6}$ &$A_{17}E_7$ & $2$ & $18$ & $[31]$\\ \hline
$N_{7}$ &$D_{10}E_7^2$ & $2$ & $18$ & $[110],[301]$\\ \hline
$N_{8}$ &$A_{15}D_9$ & $2$ & $16$ & $[21]$ \\ \hline
$N_{9}$ &$D_8^3$ & $S_3$ & $14$ & $[(122)]$ \\ \hline
$N_{10}$ &$A_{12}^2$ & $4$ & $13$ & $[15] $\\ \hline
$N_{11}$ &$A_{11}D_7E_6$ & $2$ & $12$ & $[111]$\\ \hline
$N_{12}$ &$E_6^4$ & of order $48$ & $12$ & $[1(012)]$ \\ \hline
$N_{13}$ &$A_9^2D_6$ & $2^2$  & $10$ & $[240],[501],[053]$ \\ \hline
$N_{14}$ &$D_6^4$ & $S_4$ & $10$ & $[$even perm. of $\{0,1,2,3\}]$ \\ \hline
$N_{15}$ &$A_8^3$ & $S_3\times 2$  & $9$ & $[(114)]$\\ \hline
$N_{16}$ &$A_7^2D_5^2$ & $2^3$ & $8$ & $[1112],[1721]$\\ \hline
$N_{17}$ &$A_6^4$ & $2.A_4$ & $7$ & $[1(216)]$ \\ \hline
$N_{18}$ &$A_5^4D_4$ & as $N_{12}$ & $6$ & $[2(024)0],[33001],[30302],[30033] $ \\ \hline
$N_{19}$ &$D_4^6$ & $3\times S_6$ &$6$ & $[111111],[0(02332)]$ \\ \hline
$N_{20}$ &$A_4^6$ & $2.L_2(5).2$ & $5$ & $[1(01441)]$\\ \hline
$N_{21}$ &$A_3^8$ & $2^3.L_2(7).2$ & $4$ & $[3(2001011)]$ \\ \hline
$N_{22}$ &$A_2^{12}$ & $2.M_{12}$ &$3$ & $[2(11211122212)]$\\ \hline
$N_{23}$ &$A_1^{24}$ &  $M_{24}$ & $2$ & $[1(00000101001100110101111)]$ \\ \hline
$\Lambda$ &$\emptyset$ & $Co_0$ &$0$ & $\emptyset$ \\ \hline
\end{tabular}
\end{table}

By \Ref{lem}{nik_immerge1} all of the Niemeier lattices can be defined as primitive sublattices of $\Pi_{1,25}\cong U\oplus E_8(-1)^3$ \index{Lattice, $\Pi_{1,25}$} by specifying a primitive isotropic vector $v$ and setting $N=(v^\perp\cap\Pi_{1,25})/v$.
\begin{ex}\label{ex:leech_26}
Let $\Pi_{1,25}\subset\mathbb{R}^{26}$ (the first coordinate of $\mathbb{R}^{26}$ is the positive definite one) be as before and let
\begin{align}\nonumber
v = &(17,1,1,1,1,1,1,1,1,3,3,3,3,3,3,3,3,3,5,5,5,5,5,5,5,5)\\\nonumber
w = &(70,0,1,2,3,4,5,\dots,24)
\end{align}
be two isotropic vectors in the standard basis of $\mathbb{R}^{26}$. Then 
\begin{equation}\nonumber
\Lambda\cong (w^\perp\cap\Pi_{1,25})/w
\end{equation}
 and 
\begin{equation}\nonumber 
N_{15}\cong (v^\perp\cap\Pi_{1,25})/v.
\end{equation}
\end{ex}

\subsection{The ``holy'' construction}\label{ssec:holy}

In this subsection we give a few different constructions of the Leech lattice $\Lambda$ arising from the other Niemeier lattices. These constructions are useful to classify isometries of the Leech lattice.\\ The detailed construction is contained in \cite[Section 24]{con}, in the following we just sketch it:
Let $A_n$ be a Dynkin lattice defined by
\begin{equation}\nonumber
A_n=\{(a_1,\dots,a_{n+1})\,\in\,\mathbb{Z}^{n+1},\,\,\sum{a_i}=0\}.
\end{equation}
And let $f_j$ be the vector with $-1$ in the $j-$th coordinate and $1$ in the $(j+1)-$th, zero otherwise. Let moreover $f_0=(1,0,\dots,0,-1)$. In general the $f_i$ form a set of extended roots for the Dynkin lattice.\\ Let $g_0=h^{-1}(-\frac{1}{2}n,-\frac{1}{2}n+1,\dots,\frac{1}{2}n)$, where $h$ is the Coxeter number of $A_n$ and let the $g_i$'s be a cyclic permutation of coordinates of $g_0$.
Now let $A_n(-1)^m$ be a 24 dimensional lattice and let $h_k=(g_{i_1},\dots,g_{i_m})$ where $[i_1i_2\dots i_m]$ is a glue code obtained from \Ref{tab}{nieme}. Let $f^j_i=(0,\dots,0,f_i,0,\dots,0)$ where $f_i$ belongs to the $j-$th copy of $A_n$. Let $m_i^j$ and $n_w$ be integers.\\
Then the following holds: the set of vectors satisfying
\begin{equation}\label{holy_nieme}
\sum_{j=1}^m\,\sum_{i} m_i^jf_i^j+\sum_w n_wh_w,\,\,\,\,\sum_{w}n_w=0 
\end{equation}
is isometric to the Niemeier lattice with Dynkin diagram $A_n^m$. 
While the set of vectors
\begin{equation}\label{holy_leech}
\sum_{j=1}^m\,\sum_{i} m_i^jf_i^j+\sum_w n_wh_w,\,\,\,\,\sum_{w}n_w+\sum_{i,j}m_i^j=0
\end{equation}
is isometric to the Leech lattice $\Lambda$. We call the set defined by \eqref{holy_leech} the holy construction of $\Lambda$ with hole \eqref{holy_nieme}.\\
Moreover the glue code provides several automorphisms of the Leech lattice, where the action of $t\in G(N)$ is given by sending $h_w$ to $h_{w+t}$.\\

Let us now use this construction for some isometries relevant in the following.
\begin{ex}\label{ex:p3E83}
Let us apply this construction to the lattice $E_8(-1)^3$ and let $\varphi$ be an order 3 permutation of the 3 copies of $E_8(-1)$. With the holy construction with hole $N_3$, it induces an automorphism $\varphi$ of $\Lambda$ of order 3 which fixes the only glue vector $g_0$. A direct computation shows that $T_\varphi(N_3)\cong E_8(-3)$ and $S_\varphi(N_3)\cong S_\varphi(\Lambda)=\{a-\varphi(a),\,a\in E_8(-1)^3\}$. Let us call this lattice $S_{3.exo}$ \index{Lattice, $S_{3.exo}$}, it is $S_g(\Lambda)$ for any $g\in Co_0$ in the conjugacy class $3D$ (in the notation of \cite{atlas}).
\end{ex}

\begin{ex}\label{ex:A46}
Let us apply this construction to the lattice $A_4(-1)^{6}$, we then have $G(N)=(\mathbb{Z}_{/(5)})^3$ acting on $\Lambda$. The normaliser of this group (inside $Co_1$) is one of the maximal subgroups of $Co_1$, and its structure is analysed in \cite{atlas}. The elements of $G(N)$ fall under three conjugacy classes labeled $5A,5B$ and $5C$. Each conjugacy class has respectively $40,60$ and $24$ representatives inside $G(N)$. Therefore we can compute the rank of the invariant lattice inside $\Lambda$ for each of these conjugacy classes. This rank is 4 for elements of class $5C$, 8 for elements of class $5B$ and 0 for elements of class $5A$.
\end{ex}
\begin{ex}\label{ex:A122}
Let us apply this construction to the lattice $A_{12}(-1)^{2}$, we then have $G(N)=\mathbb{Z}_{/(13)}$. Let $\varphi$ be an isometry of $\Lambda$ of order 13 generated by a non trivial element $g$ of $G(N)$ on this holy construction. Note that $\varphi$ cyclically permutes the extended roots of both copies of $A_{12}$ and therefore has no fixed points in $\Lambda$.
\end{ex}
\begin{ex}\label{ex:p11A212}
Let us apply this construction to the lattice $A_2(-1)^{12}$ and let us analyze an automorphism of order 11: it can be defined by leaving the first copy of $A_2(-1)$ fixed and by cyclically permuting the remaining 11, and the action is extended accordingly to the glue vectors. This automorphism is defined on both $N_{22}$ and $\Lambda$. Let $\varphi$ be this isometry on $A_2^{12}(-1)\otimes\mathbb{Q}$.\\
A direct computation shows that $T_{\varphi}N_{22}$ is spanned by
\begin{equation*}
f_1^1,f_1^2,\sum_2^{12}f_1^i,\sum_2^{12}f_2^i,\sum_1^{11} g_j,
\end{equation*}
where $g_j$ are generators for the glue code as in \Ref{tab}{nieme}. Keeping the same notation as before one sees that $S_{\varphi}N_{22}$ has rank 20 and is spanned by
\begin{equation}
(f^k_1-\varphi f^k_1),(f^k_2-\varphi f^k_2),(g_j-\varphi g_j).
\end{equation}

Where $k$ runs from 2 to 12. This vectors satisfy \eqref{holy_leech}, therefore this lattice is contained in $\Lambda$ and $S_{\varphi}(N_{22})=S_{\varphi}(\Lambda)$. For an explicit description of this lattice, see \cite[Example 2.9]{me_autom1}
\end{ex}

\subsection{Prime order automorphisms of the Leech lattice}\label{sec:prime_nieme}

In this subsection we give a brief analysis of prime order automorphisms on the Leech lattice, which will be used for \Ref{cor}{k3n_prime}. 
The relevant computations of this subsection can be found in \cite{me_tesi}.

\begin{prop}\label{prop:invol_nieme}
Let $\Lambda$ be the Leech lattice, $\varphi\in Aut(\Lambda)$ be an involution, then one of the following holds:
\begin{itemize}
\item $S_{\varphi}(\Lambda)=E_8(-2)$,
\item $S_{\varphi}(\Lambda)=BW_{16}(-1)$,
\item $S_{\varphi}(\Lambda))=D^+_{12}(-2)$, 
\item $S_{\varphi}(\Lambda)=\Lambda$.
\end{itemize}
Here $BW_{16}$ is the Barnes-Wall lattice, described in \cite[Section 4.10]{con}.
\begin{proof}
There are four conjugacy classes of involutions inside $Aut(\Lambda)$, with coinvariant lattices of ranks $8,12,16$ and $ 24$. The latter gives $\Lambda$, the lattice $E_8(-2)$ is obtained through the holy construction with $E_8^3$ and the lattice $D^+_{12}(-2)$ is obtained through the holy construction with $D_{12}^2$. Finally, the remaining case is given by the Barnes-Wall lattice as explained in \cite[Section 4.10]{con} 
\end{proof}
\end{prop}

\begin{prop}\label{prop:p3_nieme}
Let $\Lambda$ be the Leech lattice, and let $\varphi\in Aut(\Lambda)$ be an order three isometry, then one of the following holds:
\begin{itemize}
\item $S_{\varphi}(\Lambda)=S_{3.K3}$,
\item $S_{\varphi}(\Lambda)=W(-1)$ as in \Ref{ex}{slat_rk6},
\item $S_{\varphi}(\Lambda)=S_{3.exo}$ as in \Ref{ex}{p3E83}, 
\item $S_{\varphi}(\Lambda)=\Lambda$.
\end{itemize}
\begin{proof}
There are four conjugacy classes of order 3 isometries of the Leech lattice, whose coinvariant lattices have ranks 12,16,18 and 24 respectively. The latter is $\Lambda$, the rank 16 lattice is obtained in \Ref{ex}{p3E83} and the remaining two are obtained as orthogonal to the $\mathcal{S}$-lattice $2^{36}3^9$ and through the holy construction $E_6^4$.
\end{proof}
\end{prop}

\begin{prop}\label{prop:p5_nieme}
Let $\Lambda$ be the Leech lattice and let $\varphi\subset Aut(\Lambda)$ be an isometry of order five, then one of the following holds:
\begin{itemize}
\item $S_{\varphi}(\Lambda)=S_{5.K3}$, 
\item $S_{\varphi}(\Lambda)=S_{5.exo}$ as in \Ref{ex}{slat_5rk4}, 
\item $S_{\varphi}(\Lambda)=\Lambda$.
\end{itemize}
\begin{proof}
There are 3 conjugacy classes of order 5 elements in $Co_1$, and they can be obtained using the holy construction on $A_4^6$ as in \Ref{ex}{A46}. Keeping the same notation of that example we have that an element of class $5A$ fixes no elements of $\Lambda$. An element of class $5C$ fixes a lattice of rank 4, therefore we have $S_{5C}(\Lambda)\cong S_{5.exo}$ as in \Ref{ex}{slat_5rk4}. The last class gives $S_{5.K3}$.
\end{proof}
\end{prop}

\begin{prop}\label{prop:p7_nieme}
Let $\Lambda$ be the Leech lattice and let $\varphi\subset Aut(\Lambda)$ be an isometry of order seven, then one of the following holds:
\begin{itemize}
\item $S_{\varphi}(\Lambda)=S_{7.K3}$,
\item $S_{\varphi}(\Lambda)=\Lambda$.
\end{itemize}
\begin{proof}
There are 2 conjugacy classes $7A,7B$ of elements of order 7 and they can be both obtained by applying the holy construction to the lattice $A_6^4$ and considering automorphisms given by the glue code $G(N_{17})$: One class, such as that of the glue code $[1\,2\,1\,6]$, has $rank(S_{\varphi}(\Lambda))=24$. If we take the other class, like that of $[2\,1\,3\,0]$, we obtain $S_{7.K3}$.
\end{proof}
\end{prop}

\begin{prop}\label{prop:no13}
Let $f$ be an isometry of the Leech lattice $\Lambda$ of prime order $p$. If $p\geq 13$, $S_f(\Lambda)$ has rank at least $22$.
\begin{proof}
The only possible primes are $13$ and $23$. An automorphism of order 23 has a co-invariant lattice which is negative definite and of rank 22. This can be explicitly computed using an order 23 element of $M_{24}$ and letting it act on $\Lambda$ or on $N_{23}$. The only Niemeier lattice with an automorphism of order 13 is $\Lambda$, where all elements of order 13 are conjugate (see \cite{atlas}). These automorphisms have no fixed points on $\Lambda$, as in \Ref{ex}{A122}.

\end{proof}
\end{prop}

\end{document}